# A Comparative Analysis of Tensor Decomposition Models Using Hyper Spectral Image

Ankit Gupta [1], Ashish Oberoi [2]
Department of Computer Engineering
Maharishi Markandeshwar University
Mullana, Ambala,
Haryana - India

**ABSTRACT**

Hyper spectral imaging is a remote sensing technology, providing variety of applications such as material identification, space object identification, planetary exploitation etc. It deals with capturing continuum of images of the earth surface from different angles. Due to the multidimensional nature of the image, multi-way arrays are one of the possible solutions for analyzing hyper spectral data. This multi-way array is called tensor. Our approach deals with implementing three decomposition models LMLRA, BTD and CPD to the sample data for choosing the best decomposition of the data set. The results have proved that Block Term Decomposition (BTD) is the best tensor model for decomposing the hyper spectral image in to resultant factor matrices.

## I. INTRODUCTION

Hyperspectral image is contiguous spectrum of bands of relatively small wavelengths ranging from 0.35 to 12.5 micrometer. It measures spectral reflectance of at a series of narrow bands. These images are used to distinguish materials which are spectrally similar. Most of the hyperspectral sensors are airborne with certain exceptions. Due to six types of gaseous disturbances present in space, rays received at the sensor can be less or more than the light reflected from the surface. Various atmospheric correction algorithms have been designed to measure the concentration of gases from the hyperspectral data. Various methods like whole pixel, Sub-pixel methods, matched filtering etc. have been developed for extracting detailed information from hyperspectral imagery [11].

Data compressive representation, spectral signature identification of ingredient materials and their fractional abundances determination are the three most important objectives that should be attained for processing hyperspectral data related to an object. Applications of hyperspectral imaging (HSI) includes identification of minerals, buildings and other objects on earth, space object identification[7], face detection[13], planetary exploitation [14], biology[15]. Analyzing hyperspectral data is quite challenging because of high noise content and spatial and spectral deblurring. HSI presents substantial amount of spectral redundancy by the producing signature corresponding to each pixel in an image among many highly correlated bands. Dimensionality reduction is one of the main issues for image analysts since it is necessary to remove noise from the original data. [8].

The main problem for target detection, image segmentation and pixel classification are the pixels located at high dimensional space which enhances the computational complexity and degrade accuracy whereas spatial resolution must be good for identification of relatively small constituent materials. Otherwise it may lead to noisy image map if the spatial resolution is not good [8].Some of the hyperspectral sensors that are capturing a continuum of images of the earth surface from different angles are AVIRIS (Airborne Visible-Infrared Imaging Spectrometer) by NASA covering wavelength range from 0.4 to 2.5 micrometer using nominal spectral resolution of 10nm and operating on 224 spectral channels.

These sensors capture GBs of data every day. Vast number of techniques has been used for the analysis of hyperspectral data. These techniques are intrinsically full pixel or mixed pixel technique. End members can be identified using full pixel techniques whereas mixed pixel techniques deals with identification of pixels describing more than one element. For mixed pixel techniques spectral unmixing is one of efficient mechanism which deals with decomposition of mixed pixels to end members or pure spectra [9]. Hence identification of pure spectra or pure pixels is one of the major objectives for analysis. Due to multidimensional nature of hyperspectral image collected by airborne/satellite sensors multi-arrays can be used for its analysis. Multi-way arrays are also called tensor. The term tensor was coined by William Rowan Hamilton in 1846[5].Tensor is also called multidimensional vectors and can be used for modeling and interpretation of multidimensional data. Applications of tensor include hyperspectral image processing [6] [8], signal processing [19], neuroscience [23], graph analysis [21], computer vision [22], data mining [20] etc. [5]. Simplified mathematical approach to define a tensor as a mathematical object described by finite set of indices satisfying multi-linearity. These indices are called tensor orders.

## II. TENSOR DECOMPOSITION

Tensor decompositions were invented by Hitchcocki in 1927[16, 17].It helps in factorizing a tensor in to factor matrices in a way that resulting matrices could provide mixed data in physically realizable representation [6]. This paper deals with implementation of three tensor decomposition models namely Canonical Polyadic Decomposition (CPD), Block Term Decomposition, Low multi-linear rank approximation (LMLRA) using a hyperspectral image.

*A. Notations*

A real valued tensor defined by $\mathbb{R} \in K^{I_1 \times I_2 \times I_3 \ldots \times I_n}$ having elements or scalar components $r_{i_1 i_2 i_3 \ldots i_n}$ with $i_{j=1, 2, \ldots, I_j}$. Index $I_j$ is called a mode of tensor. Third order tensor will have 3 modes $I_1 \times I_2 \times I_3$. Sub-tensors defined by $R(l_1, l_2, \ldots, l_N)$ of orignal tensor A obtained by restricting indices to belong to subsets $K_n \subseteq \{1, 2, \ldots I_n\}$ having vector values are called fibers and matrix valued sub tensors are called slices. Fibers are created by fixing all indices except one index where as slices are acted by fixing all except two indices. Mathematically $r(:, i_2, i_3, \ldots, i_N)$ defines a fiber and $R(:,:, i_3, \ldots, i_N)$ specifies the tensor slice obtained by fixing some indices of a tensor R.Tensor manipulation deals with reshaping or reformatting of tensor. It is also called tensor folding or matricization for BTD tensor model. Multiplication of a tensor A of mode n with a matrix B is the multiplication of all mode-n vector fibers with matrix B. The diagonal tensor is defined by $Z = diag_N(\lambda_1, \lambda_2, \lambda_3, \ldots, \lambda_R)$ of order N with $d_{rr\ldots r} = \lambda_r$.

*B. Canonical polyadic decomposition (CPD)*

The Canonical Polyadic Decomposition (CPD) provides an approximation of a higher order tensor with a sum of *R* rank-one tensors. Mathematically, for a third order tensor defined by $\mathbb{R} \in K^{I_1 \times I_2 \times I_3}$ CPD is the linear combination of rank-one tensors in the form

$$\mathbb{R} = \sum_{r=1}^{R} \lambda_r b_r^{(1)} \circ b_r^{(2)} \circ b_r^{(3)} \circ b_r^{(4)} \ldots \circ b_r^{(N)} \quad (1)$$

Equivalently, $\mathbb{R}$ can be expressed as a multi-linear product with a diagonal core:

$$\mathbb{R} = Z X_1 B^{(1)} X_2 B^{(2)} X_3 B^{(3)} X \ldots X_N B^{(N)}$$
$$= [\![Z; B^{(1)}, B^{(2)}, B^{(3)}, \ldots, B^{(N)}]\!] \quad (2)$$

where $Z = diag_N(\lambda_1, \lambda_2, \lambda_3 \ldots \lambda_R)$. The smallest value of R for which (3) holds exactly will be the tensor rank and CPD is the minimum rank polyadic decomposition. CPD is also known as CANDECOMP/PARAFAC decomposition. For canonical polyadic decomposition, matrix/vector form can be obtained via the Khatri-Rao products as:

$$\mathbb{R}_{(n)} = B^{(n)} D (B^{(n)} \odot \ldots \odot B^{(n+1)} \odot B^{(n-1)} \odot \ldots \odot B^{(1)})^T \quad (3)$$

$$vec(\mathbb{R}) = [B^{(n)} \odot \ldots \odot B^{(n+1)} \odot B^{(n-1)} \odot \ldots \odot B^{(1)}] Z \quad (4)$$

Where $Z = (\lambda_1, \lambda_2, \lambda_3 \ldots \lambda_R)^T$. For the representation of higher order tensor, the number of complex-valued rank-1 terms can be strictly less than the number of real-valued rank-1 terms [8] the determination of tensor rank is in general NP-hard. CORCONDIA algorithm [18] is one of the existing techniques for rank estimation based on diagonality of a core tensor. However uniqueness conditions play an important role for exact tensor decomposition by providing it theoretical bounds. Uniqueness condition for a third order tensor [10], states that CPD is unique up to unavoidable scaling and permutation ambiguities, provided that

$$k_{B^{(1)}} + k_{B^{(2)}} + k_{B^{(3)}} \geq 2R + 2 \quad (5)$$

where the Kruskal rank $k_B$ of a matrix B is the maximum value ensuring the linear dependency of any subset of $k_B$ columns. Various tensor unfolding models like Tucker Decomposition model, Canonical polyadic Decomposition (CPD) etc. are available to process hyperspectral data. In a case when factor matrix has full column rank, more relaxed uniqueness conditions can be obtained. Under more natural and relaxed conditions that require the components to be "sufficiently different", CPD would be unique.

### C. *Low multi-linear rank approximation*

LMLRA decomposes large-scale data tensor $\mathbb{R}$ in to a core tensor multiplied by set of factor matrices along each mode in an approximate manner [12].It is similar to high order principal component analysis. For a third order tensor, LMLRA is defined by

$$\mathbb{R} \approx E X_1 A X_2 B X_3 C = \sum_{i=1}^{I}\sum_{j=1}^{J}\sum_{w=1}^{W} e_{ijw} a_i \circ b_j \circ c_k = [\![E:A,B,C]\!] \quad (6)$$

The factor matrices A, B, C are the principal components among in each mode and E are the core tensor; i, j, w are the components or columns in the factor matrices. Matricized form of tensors one in each mode are:

$$R_1 = A E_{(1)} (C \otimes B)^T \quad (7)$$

$$R_2 = B E_{(1)} (C \otimes A)^T \quad (8)$$

$$R_3 = C E_{(1)} (B \otimes A)^T \quad (9)$$

### D. *Block term decompositions*

It is a combination of tucker and CP. In this model tensor is represented as sum of low rank tucker tensors. Replacement of rank-1 matrices $b_r^{(1)} \circ b_r^{(2)} = b_r^{(1)} b_r^{(2)T}$ by low rank matrices $A_s B_s^T$ in (3) leads to the representation of tensor $\mathbb{R}$ as

$$R = \sum_{s=1}^{S} (A_s B_s^T) \circ C_s \quad (10)$$

If we use terms required for having low multi-linear rank then the equations become

$$R = \sum_{s=1}^{S} G_s X_1 A_s X_2 B_s X_3 C_s \quad (11)$$

This is called as Block Term Decomposition. It confesses the modeling of more complex signal components as compared to CPD and is more unique under more restrictive but fairly natural conditions.

### III. RELATED WORK

The key aspect of processing an image is to define and extract the information from it. There are two ways of information definition: 1) Recognizing objects in an image, 2) deriving an entity which represents the information in an image **[2].** Processing hyperspectral data deals with the estimation of missing values in the three dimensional matrix. The basic problem in missing value estimation is to find the relationship between known and unknown elements **[3]**.Due to the multidimensionality of hyperspectral data multi-way data analysis can be used to extract the hidden structures and correlate the elements within a multi-way array. These multi-way arrays are often called tensors. The term tensor is used as high order generalization of vectors and matrices. Each dimension of a tensor is called mode and number of variables in each dimensions indicates its dimensionality. A tensor data model consists of two parts: 1) Structural part 2) Residual part. Fitness of a model depends on the analysis of its residual part. Most of the tensor unfolding models preserves the multi-way nature of data e.g. Tucker decomposition which deals with splitting a tensor in to matrices. This process may lead to loss of information and misinterpretation of highly noised data. [4] The other model called CPD, Canonical Polyadic Decomposition historically known as CANDECOM/PARAFAC model which deals with decomposing a tensor using R rank-one terms [5].These models can be applied to any multivariate data for extracting useful information [4]. Methods like High Order Partial Least Square (HOPLS) can be used to predict values of data from real scene for plotting it as an image map virtually in low dimensional common latent subspace [1].

### IV. RESULTS

This section deals with exploring the efficacy of the tensor models explained in section 2 for decomposition of third order tensor in to its factors specified by particular model. It deals with selecting the best tensor decomposition model for our data.The performance of the decomposition models was evaluated on the basis residual error produced during decomposition and relative error generated depending on the type of model. The algorithms performed number of iterations in order to minimize the resulting residual error whereas relative error is the relative error generated with refinement step of algorithms for best decomposition of the input data.

### A. *CPD*

CPD decomposes third order tensor to rank-1 tensors. The algorithm is slow because computing CPD includes data compression, generating random factor matrices and refinement using CPD alternating least square method. The variation of residual error corresponding to CPD is shown in figure 1. Relative error generated using compression, pseudo random matrix generation are 0.129035, 1.00063 respectively. CPD computation and refinement uses alternating least square method and errors generated during the steps are 0.0606507, 0.0828356 respectively.

### B. *LMLRA*

LMLRA decomposes a tensor in to core tensor whose size is less than original tensor and factor matrices The algorithm consist of two steps 1) generating matrices 2) LMLRA computation using non-linear least square .The relative error corresponding to the generation of matrices and LMLRA computation is 1.80507, 0.0454302 respectively .The residual error graph corresponding to method is shown in figure3.

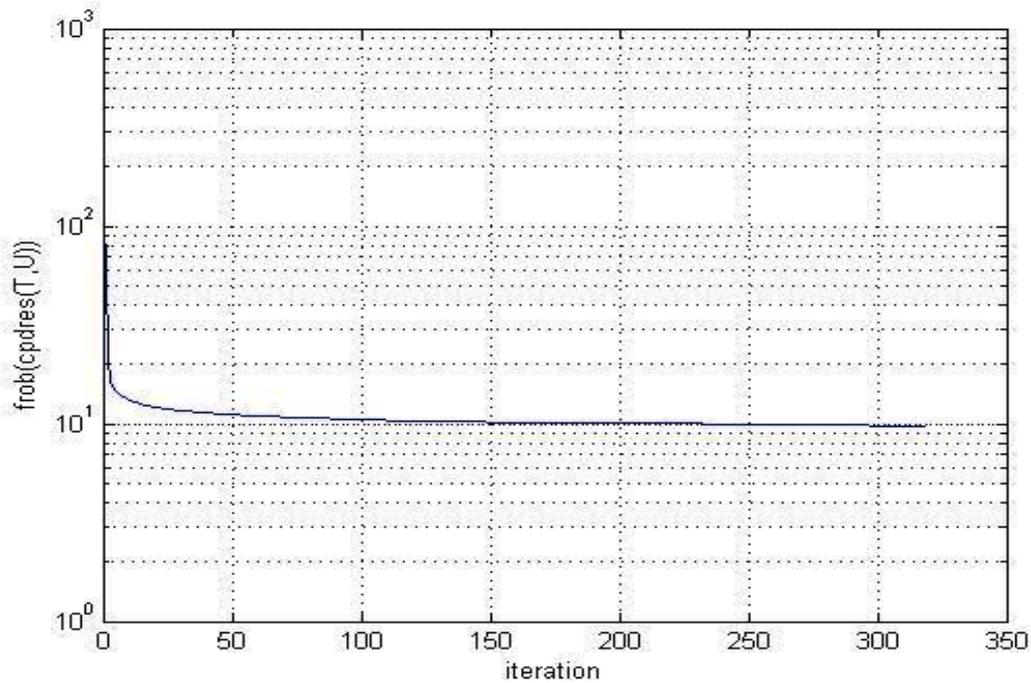

Fig 1: Residual during Decomposition graph generated during CPD computation.

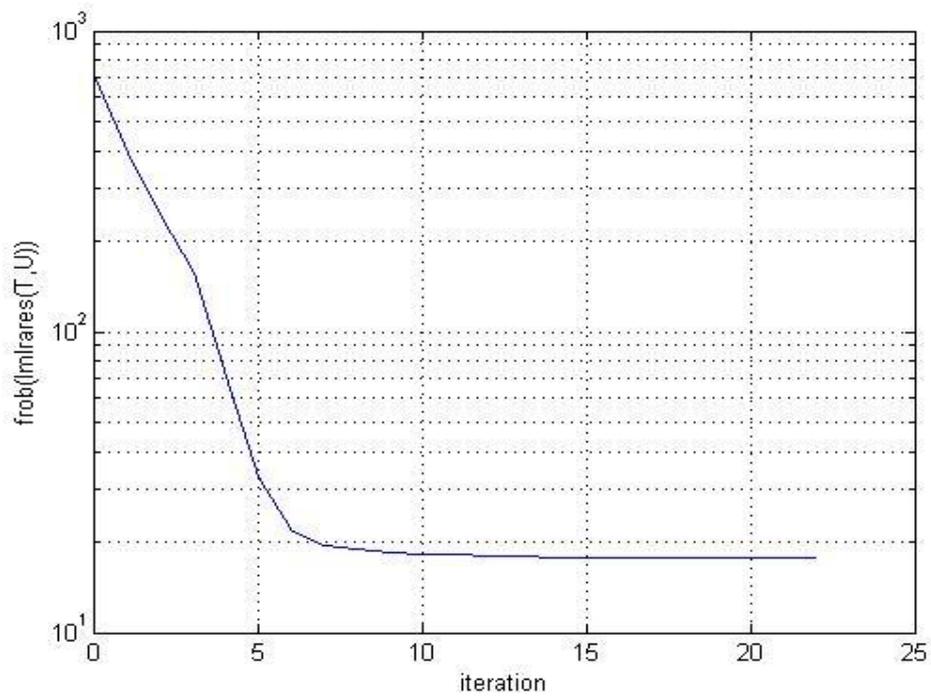

Fig 2: Residual during Decomposition graph generated during LMLRA computation.

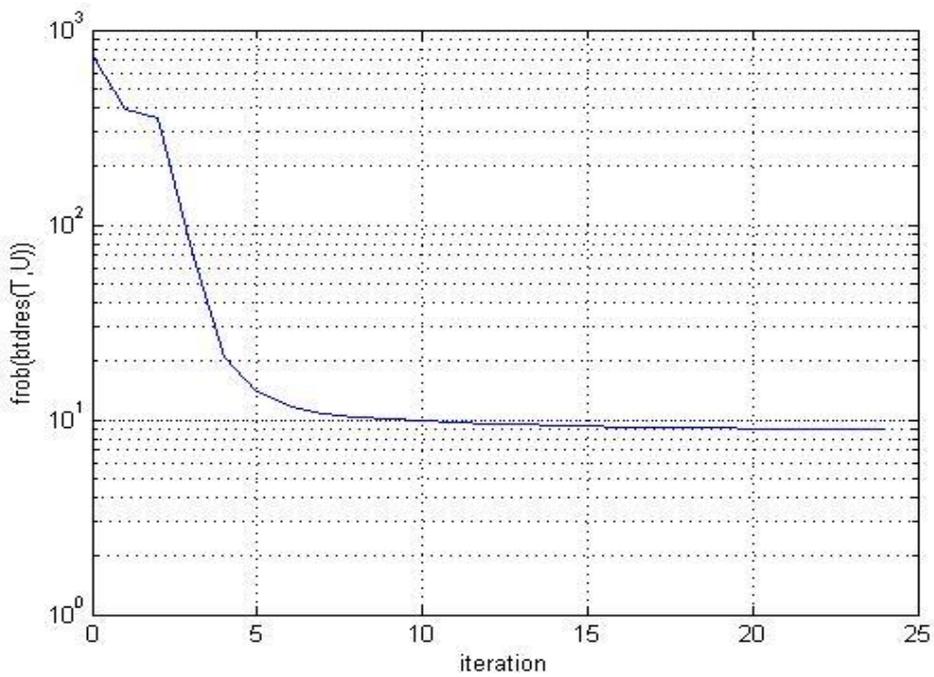

Fig 3: Residual during Decomposition graph generated during BTD computation.

C. *BTD*

Block Term Decomposition is the combination of Canonical Polyadic Decomosition and tucker decomposition.This method was proved to be fast as compared to other two decomposition models.The algorithm used for BTD computation is non-linear least square method .The relative error corresponding to this method is 0.022868829769866. The residual error graph for BTD is shown in figure 2.

TABLE I
NUMBER of ITERATIONS and RELATIVE ERROR for DECOMPOSITION MODELS.

| Methods | iterations | Relative Error |
|---------|-----------|----------------|
| CPD     | 318       | 0.0828355865100 |
| BTD     | 24        | 0.0228688297698 |
| LMLRA   | 22        | 0.0455030082908 |

## V. CONCLUSION

The Relative errors of CPD, BTD, and LMLRA are shown in the table 1. Based on the comparison table and residual error graphs corresponding to each model BTD has been proved to be the best tensor decomposition model for our data set since relative error of BTD is the smallest among all models. On the basis of upper bound of residual error CPD is best where BTD performed the best on the basis of lower bound of residual error. CPD took long time to adjust its residual error to minimum as number of iterations performed is 318 whereas the iterations performed by LMLRA and BTD 22, 4 respectively are enormously less than canonical polyadic decomposition.

It is concluded that BTD was selected as the best decomposition model for our hyperspectral image cube based on lower value relative error and final residual (lower bound) whereas CPD gave the worst performance against our data set beside the lowest value of upper bound of residual than other models. Future work includes creating High order partial least square method which could take block term decomposition outputs as the input for predicting values from the original sample thereby plotting them virtually on to a latent subspace, implementing multi-linear PCA on the data set for further hyperspectral image processing.